\input amstex
\magnification=\magstep1
\input amsppt.sty
\topmatter
\title{Approximate factorization and
concentration for characters of symmetric groups}\endtitle
\author Philippe Biane 
\endauthor
\leftheadtext {Philippe Biane }\rightheadtext{Approximate factorization and
concentration}
\address CNRS,
DMA, \' Ecole Normale
Sup\'erieure\endgraf
45, rue d'Ulm, 75005 Paris FRANCE\endgraf
{\rm Philippe.Biane\@ens.fr}
\endaddress
\abstract We prove a factorization-concentration result for characters of
symmetric groups. This is then applied to the asymptotic behaviour of the
decomposition of the tensor representations. There are 
 connections
with the Pastur-Marcenko distribution of random matrix theory,
 and freely infinitely divisible
distributions.
\endabstract
\thanks
The author thanks S. Kerov and A. Okounkov
 for several enlightening discussions.
\endthanks
\endtopmatter
\document
\head{ Introduction}\endhead
 In \cite{B} we
determined the asymptotic behaviour of normalized irreducible 
characters of the symmetric
groups $S_q$, as $q\to\infty$, under the hypothesis that the corresponding Young
diagrams, rescaled by a factor $q^{-1/2}$, have a limit shape. It turns out that
for such characters $\chi$, the order of magnitude of $\chi(\rho)$ is
$q^{-\vert\rho\vert/2}$ where the length $|\rho|$
 of a permutation $\rho$ is the minimal
number $n$ required to write  $\rho$ as a product of $n$ transpositions.
Alternatively, $|\rho|=q-c(\rho)$ where $c(\rho)$ is the number of cycles of 
$\rho$. Furthermore these characters satisfy an approximate factorization
property
$$\chi(\rho\sigma)=\chi(\rho)\chi(\sigma)+o(q^{-(\vert
\rho\vert+\vert\sigma\vert)/2})\tag 0.1$$
for permutations $\rho,\sigma$ with disjoint supports.
In this paper we consider normalized positive
definite functions
 on symmetric groups. Such a function defines a probability measure on the set
 of Young diagrams. We shall consider functions $\psi$ satisfying bounds
$$\vert\psi(\rho)\vert\leq c_nq^{-n/2}\tag 0.2$$
for $\vert\rho\vert=n$. We prove that for large $q$, 
 the  probability measure associated with a 
 $\psi$ satisfying the approximate factorization (0.1) is concentrated
  on the set of Young diagrams
whose normalized character takes values close to that of $\psi$.
 As an application of this result we shall
 investigate the asymptotic behaviour of the decomposition of the tensor
 representation $(\Bbb C^N)^{\otimes q}$ of the symmetric group $S_q$. We shall
 prove that in the asymptotic regime $q\to\infty$,
 ${ \sqrt q\over N}\to c\in [0,+\infty[$
  the typical Young diagram occuring in the
 decomposition of this representation has a certain limit shape depending on
 $c$. For $c=0$ we recover the limit shape occuring in the asymptotics
 of the Plancherel measure on $S_q$,
  as in the results of Kerov-Vershik \cite{KV} and
 Logan-Shepp 
 \cite{LS}.  For $c>0$ this limit shape 
 is strongly related to the
 so-called Pastur-Marcenko distribution occuring in the theory of random
 matrices, and in the theory of free probability, where it appears under the
 name of ``free Poisson distribution''. We shall also consider more general
 measures on Young diagrams coming from tensor states on 
 $\Cal B(H)^{\otimes q}$. By tuning the state with $q$
  we will get limit diagrams which are related to  freely infinitely
  divisible
  measures.\par
  This paper is organized as follows. In
 Section 1 we settle notations, recall results from \cite{B},
  and state the main result of the paper. In
 Section 2 we prove the result, and in 
 Section 3 we discuss the decomposition of
 the tensor representation.
 \head{1. Notations and statement of the main result}\endhead 
 We recall notations from \cite{B}, to which we refer for more details.
 \subhead 1.1 Symmetric groups
and  Young diagrams\endsubhead
 We denote by $S_q$
 the symmetric group on $\{1,\ldots,q\}$,  
  by $(ij)$ the transposition exchanging $i$ and $j$, and more
generally by $(i_1i_2\ldots i_n)$  cyclic permutations of order $n$. 
For $\sigma\in
S_q$, let $c(\sigma)$ be its number of 
 cycles, and $s(\sigma)$ be the number of elements not fixed by $\sigma$. Then
$\vert\sigma\vert:=q-c(\sigma)$ is the smallest number $n$ such that
$\sigma$ can be written as a product of $n$ transpositions. 
 We denote by $\Cal Y_q$ the set of Young diagrams with $q$ boxes, and
$\Cal Y=\cup_{q=1}^{\infty}\Cal Y_q$.
  If $\lambda\in \Cal Y_q$ let
$[\lambda]$ be the associated 
irreducible representation of $S_q$, and $\chi_{\lambda}$ be its normalized
character, i.e. $\chi_{\lambda}(e)=1$.  
Recall from \cite{B}, \cite{K} that a Young diagram can be identified with
 a piecewise linear function with slopes $\pm 1$,  and  local
minima
and maxima occuring at two interlacing sequences of integer points
$$x_1<y_1<x_2<\ldots<y_{n-1}<x_n$$
as in the following example
\input epsf.tex
$$\aligned{\epsfbox{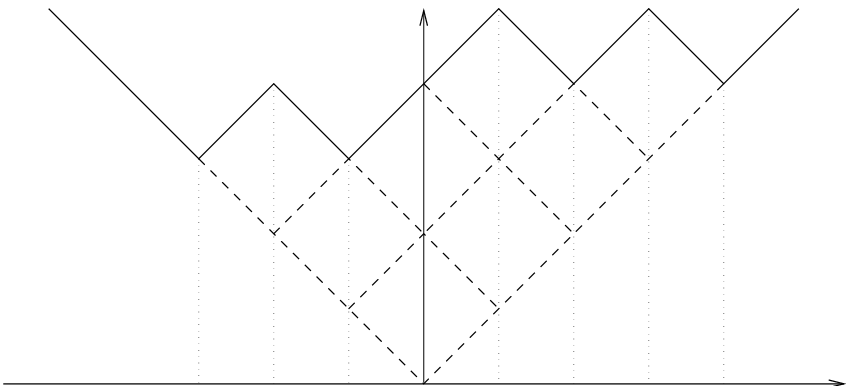}}\\ x_1\quad\,  y_1\quad\, 
x_2\qquad\quad\,\,  y_2\quad\, x_3\quad\, y_3\quad\, x_4\qquad\quad\endaligned$$
and that we can
 embed the set of Young diagrams in the   
 space $\Cal C\Cal Y$ of continuous Young diagrams, i.e. 
  functions $\omega:\Bbb R\to
\Bbb R$ satisfying 
\roster\item $\vert\omega (u_1)-\omega(u_2)\vert\leq \vert
u_1-u_2\vert\quad \text{for all}\quad u_1,u_2\in\Bbb R$
\item $\omega(u)=\vert u\vert$ for sufficiently large $\vert
u\vert$.\endroster
To each continuous diagram $\omega\in\Cal C\Cal Y$ one can associate
a probability measure $\frak
m_{\omega}$ with compact support on $\Bbb R$,  determined by the equation
$${1\over z}\exp\int _{\Bbb R}{1\over
x-z}\sigma'(x)\,dx= \int_{\Bbb R}{1\over z-x}\frak
m_{\omega}(dx)\quad\quad z\in\Bbb C\setminus \Bbb R\tag 1.1$$
where
$\sigma(u)=(\omega(u)-\vert u\vert)/2; u\in \Bbb R$. 
This measure is 
called the transition measure of the diagram and its moments are
called the moments of the corresponding continuous diagram. 
We shall denote them
by $$m_n(\omega)=\int_{\Bbb R} x^n\,\frak m_{\omega}(dx)$$
The measure corresponding to a diagram in $\Cal C\Cal Y$ has bounded support, 
it   is centered (i.e.
$m_1=0$), and its 
second moment is equal to half the area of the set $$\{(u,v)\in\Bbb
R^2\vert\,
|u|\leq v\leq \omega(u)\}$$
Besides moments of measures we shall also consider their free cumulants. Recall
 (\cite{B}, \cite{VDN}),
 that the free cumulants are defined as the coefficients $R_n$ in the expansion
$$K(z)={1\over z}+\sum_{n=1}^{\infty}R_nz^{n-1}$$
where $K(z)$ is the 
functional 
inverse of the Cauchy transform $$G(z)=\int_{\Bbb R}{1\over z-x}\frak
m_{\omega}(dx)$$
\subhead{1.2 Positive definite functions}\endsubhead
Let $\psi$ be a normalized (i.e. $\psi(e)=1$), central
 positive definite function on
$S_q$.
By Fourier analysis on $S_q$, one can
 expand $\psi$ as a convex linear combination of normalized
characters
$$\psi=\sum_{\lambda\in\Cal Y_q}p_{\lambda}\chi_{\lambda}$$
The weights $p_{\lambda}$ are non negative and sum to one, 
therefore they define a probability measure
$\Pi_{\psi}$ on $\Cal Y_q$, which puts a mass $p_{\lambda}$ on $\lambda$.
 Let $\Cal V\subset \Cal Y_q$ and $\gamma>0$,
 we say that $\psi$ is {\it $\gamma$-supported } on $\Cal V$ if $\Pi_{\psi}(\Cal
 V)>1-\gamma$.
 
 Let us define a probability measure on $\Bbb R$
  associated with $\psi$ by
 $$\frak m_{\psi}=\sum_{\lambda\in\Cal Y_q}p_{\lambda}\frak m_{\lambda}$$
 and denote  its moments by 
 $$m_n(\psi)=\sum_{\lambda\in\Cal Y_q}p_{\lambda}m_n(\lambda)$$
 To this probability measure we can associate a diagram, by the correspondence
 (1.1), which we call
 $\omega_{\psi}$.
By the GNS construction,  
 there exists a finite dimensional complex Hilbert  space $V$, a unitary
 representation $r_{\psi}:S_q\to \Cal B(V)$ (where $\Cal B(V)$ is the space of
 bounded operators on $V$), and a state $\tau_{\psi}$ on $\Cal B(V)$, tracial
  on $r_{\psi}(S_q)''$
 such that 
 $$\psi(\rho)=\tau_{\psi}(r_{\psi}(\rho))$$
for all $\rho\in S_q$.   
As in \cite{B} for a unitary representation $r$ of $S_q$ in some $\Cal B(H)$,
we define 
the selfadjoint element $$\Gamma(r)=\pmatrix 0&1&1&1&\ldots &1&1\\
1&0&r(1\,2)&r(1\,3)&\ldots &r(1\,q-1)&r(1\,q)\\
1&r(1\,2)&0&r(2\,3)&\ldots &r(2\,q-1)&r(2\,q)\\
\vdots&\vdots&\vdots&\vdots&\vdots&\vdots&\vdots\\
1&r(1\,q)&r(2\,q)&r(3\,q)&\ldots&r(q-1\,q
)&0\endpmatrix$$
of $\Cal B(H)\otimes  M_{q+1}(\Bbb C)$. By \cite {B}, Proposition 3.3, 
the measure $\frak m_{\psi}$ is
the distribution of the self adjoint element $\Gamma(\psi):=\Gamma (r_{\psi})$
in the non-commutative probability space 
$$(\Cal B(V_{\psi})\otimes M_{q+1}(\Bbb C), \tau_{\psi}
\otimes \langle.\rangle)$$ 
where $\langle.\rangle$ denotes the normalized trace on 
$M_{q+1}(\Bbb C)$. 
In particular, one has
$$m_n(\psi)=\tau_{\psi}(\langle \Gamma(\psi)^n\rangle)\qquad\text{for all $n\geq
1$}$$
\subhead{1.3 Approximate factorization}\endsubhead
We now give a precise definition of the approximate factorization property
(0.1).
\proclaim{Definition 1.3}
Let $n$ be a
 positive integer,
 let $c=(c_l)_{l\geq 1}$ be a sequence of positive real numbers,
  and  let $\delta>0$, 
  define $\Cal F_{c,\delta}^n$ as the set of all positive definite normalized
 functions $\psi: S_q\to \Bbb C$, for some positive
  integer $q$, such that 
 \roster 
 \item $\vert \psi(\rho)\vert\leq c_lq^{-l/2}$ for all $l\leq n$, and $\rho\in
 S_q$
 with $\vert\rho\vert= l$.
 \item $\vert \psi(\rho\sigma)-\psi(\rho)\psi(\sigma)\vert\leq \delta
 q^{-l/2}$ for all $l\leq n$ and all $\rho,\sigma\in S_q$, with disjoint
 supports, such that $|\rho\sigma|=l$.
 \endroster
 \endproclaim
 Note that $|\rho\sigma|=|\rho|+|\sigma|$ if $\rho$ and $\sigma$ have disjoint
 supports, so that (2) means that, according to (1),
  $\psi(\rho\sigma)$ has the same order of magnitude as
  $\psi(\rho)\psi(\sigma)$.
 The main result of this paper is the following
 \proclaim{Theorem 1} Let the sequence $(c_l)_{l\geq 1}$ and the integer
 $n$ be given, then there exists a constant
 $K$, such that for all $\delta\geq 0$,  every
  $\psi\in\Cal F_{c,\delta}^n$ is $(\delta+q^{-1})^{1/2}$-supported by
   the set of Young diagrams $\lambda\in\Cal Y_q$
  satisfying 
  $$|m_l(\lambda)-m_l(\psi)|\leq
 K(\delta+q^{-1})^{1/4}q^{l/ 2} $$ for all
 $l\leq n/2$.
 \endproclaim
 Let us denote by $R_n(\psi)$ the free 
 cumulants of the measure $m_{\psi}$, then using the moment-cumulant formula
 of \cite{S} (see also Section 2 of \cite{B}), we
 get 
 the following corollary of the previous result
 \proclaim{Corollary 1 }
  Let the sequence $(c_l)_{l\geq 1}$ and the integer
   $n$ be given, then there exists a constant
 $K$, such that for all $\delta\geq 0$,  every
  $\psi\in\Cal F_{c,\delta}^n$ is $(\delta+q^{-1})^{1/2}$-supported by
   the set of Young diagrams $\lambda\in\Cal Y_q$
  satisfying 
  $$|R_l(\lambda)-R_l(\psi)|\leq
 K(\delta+q^{-1})^{1/4}q^{l/ 2} $$ for all
 $l\leq n/2$.
 \endproclaim
 Using the results of \cite{B} and the above Corollary, 
 a reformulation of Theorem 1 can be given which
 involves character values instead of moments.
 \proclaim{Theorem 2} Let the sequence $(c_l)_{l\geq 1}$ and the integer
 $n$ be given, then there exists a constant
 $K$, such that for all $\delta\geq 0$, every
  $\psi\in\Cal F_{c,\delta}^n$ is $(\delta+q^{-1})^{1/2}$-supported 
  by the set of Young diagrams $\lambda\in\Cal Y_q$
  satisfying 
$$|\chi_{\lambda}(\rho)-\psi(\rho)|\leq
 K(\delta+q^{-1})^{1/4}q^{-l/2} $$ for all $\rho\in S_q$ with 
 $|\rho|=l\leq n/2$.
 \endproclaim

 \head {2. Proof of Theorems 1 and 2}\endhead
 We shall rely heavily on results of \cite{B}, especially Section 4.
  Let $\psi$ be a normalized
 positive definite function, then $m_n(\psi)=\sum_{\lambda\in \Cal Y_q}
 p_{\lambda}m_n(\lambda)$ is the mean value for $\Pi_{\psi}$
  of the random variable $\lambda\mapsto m_n(\lambda)$ on $\Cal Y_q$. The proof
  of Theorem 1 consists
   in giving an estimate for the variance of this random
  variable and using Markov's inequality. 
 The key result is the following whose proof is analogous to the proof of Lemma
 5.1.1 or 6.1 in \cite{B}.
 \proclaim{Lemma 2.1}
 Let $n_1,n_2,\ldots n_k$ be positive integers, then for every sequence
 $(c_l)_{l\geq 1}$ and $\delta>0$
 there exists a constant $K$ such that
 for all $\psi\in \Cal F_{c,\delta}^{n}$ with $n=n_1+\ldots +n_k$, one has
 $$\gather\left|\tau_{\psi}(\langle \Gamma(\psi)^{n_1}\rangle
 \langle \Gamma(\psi)^{n_2}\rangle\ldots 
 \langle \Gamma(\psi)^{n_k}\rangle)-
 \tau_{\psi}(\langle \Gamma(\psi)^{n_1}\rangle)
 \tau_{\psi}(\langle \Gamma(\psi)^{n_2}\rangle)\ldots 
 \tau_{\psi}(\langle \Gamma(\psi)^{n_k}\rangle)\right|\\
 \leq K(\delta+q^{-1})q^{n/2}
 \endgather $$
 \endproclaim
\demo{Proof}
We shall only need the case $k=2$ of this Lemma which is what we prove here.
The argument can be easily extended to yield the general case, 
which can be used to give estimates on higher moments of
$\langle\Gamma(\psi)^n\rangle-\tau_{\psi}(\langle\Gamma(\psi)^n\rangle)$.
\par
 One
has
$$
\langle\Gamma(\psi)^n\rangle={1\over q+1}\sum_{0\leq i_1,i_2,\ldots ,i_n\leq
 q}r_{\psi}((i_1i_2)\ldots (i_{n-1}i_n)(i_ni_1))
$$
and analogously
$$
{\aligned 
&\langle\Gamma(\psi)^{n_1}\rangle
\langle\Gamma(\psi)^{n_2}\rangle=\\
{1\over (q+1)^2}&
\sum_{0\leq i_1,i_2,\ldots ,i_{n_1}\leq
 q\atop 0\leq j_1,j_2,\ldots ,j_{n_2}\leq
 q}r_{\psi}\bigl((i_1i_2)\ldots(i_{n_1}i_1)
 (j_1j_2)\ldots(j_{n_2}j_1)\bigr)
 \endaligned}
 \tag 2.1.1
 $$
 where by convention $(ii)$ is zero and $(ij)=1$ if either $i$ or $j$ (but not
 both) is 0.
 Applying the state $\tau_{\psi}$ one gets
$$
{\aligned 
&\tau_{\psi}(\langle\Gamma(\psi)^{n_1}\rangle
\langle\Gamma(\psi)^{n_2}\rangle)=\\
{1\over (q+1)^2}&
\sum_{0\leq i_1,i_2,\ldots ,i_{n_1}\leq
 q\atop 0\leq j_1,j_2,\ldots ,j_{n_2}\leq
 q}\psi\bigl((i_1i_2)\ldots(i_{n_1}i_1)
 (j_1j_2)\ldots(j_{n_2}j_1)\bigr)
 \endaligned}
 \tag 2.1.2
 $$
 We shall compare (2.1.2) to 
$$
{\aligned
&\tau_{\psi}(\langle\Gamma(\psi)^{n_1}\rangle)\tau_{\psi}(
\langle\Gamma(\psi)^{n_2}\rangle)=\\
{1\over (q+1)^2}&
\sum_{0\leq i_1,i_2,\ldots ,i_{n_1}\leq
 q\atop 0\leq j_1,j_2,\ldots ,j_{n_2}\leq q}
 \psi\bigl((i_1i_2)\ldots (i_{n_1}i_1)\bigr)\,
 \psi\bigl( (j_1j_2)\ldots(j_{n_2}j_1)\bigr)
 \endaligned}
 \tag 2.1.3
 $$ 
 In (2.1.2), since $\psi$ is a central function, the value of
 $$
 \psi\bigl((i_1i_2)\ldots (i_{n_1}i_1)
 (j_1j_2)\ldots (j_{n_2}j_1)\bigr)
 $$ depends only on the set $J$ of places where $i_k$ or $j_k$ is 0, and the
 partition of the set 
 $(\{1,2,\ldots,n_1\}\cup\{1',2',\ldots,n_2'\})\setminus J$ given by
 the equivalence relation $u\sim v$ if $i_u=i_v$, $u\sim v'$ if $i_u=j_{v'}$, or
 $u'\sim v'$ if $j_{u'}= j_{v'}$. Given such a set $J\subset
 \{1,2,\ldots,n_1\}\cup\{1',2',\ldots,n_2'\}$ and the  partition $\pi$,
  the number of corresponding terms
 in the sum is equal to $(q)_{comp(\pi)}$ where $(q)_r=q(q-1)\ldots (q-r+1)$
  and
 $comp(\pi)$ is the number of components of $\pi$. Denoting by
 $\psi(\pi,J)$ the common value of $\psi$ on sequences corresponding to the set
 $J$ and the 
 partition $\pi$,   the right hand side of (2.1.2) can be rewritten as 
  a sum  of the form
 $$
 \sum_{\pi,J}{1\over (q+1)^2}(q)_{comp(\pi)}\psi(\pi,J)\tag 2.1.4
 $$
 Analogously we can write the right hand side of (2.1.3) as
  $$
 \sum_{\pi_1,J_1\atop\pi_2,J_2}
 {1\over (q+1)^2}(q)_{comp(\pi_1)}\psi(\pi_1,J_1)
 (q)_{comp(\pi_2)}\psi(\pi_2,J_2)\tag 2.1.5
 $$
 Let $h(\pi)$ denote the conjugacy class of the permutation corresponding to
 $\pi$, and let $|h(\pi)|$ denote the length of any of its elements. 
 Hypothesis (1) of Definition 1.3 implies a bound   
 $|\psi(\pi,J)|\leq C q^{-|h(\pi)|/2}$. Furthermore, by an argument similar to
 Lemma 4.3.2 of \cite{B} one can see that if $J\not=\emptyset$, then 
 $$n+|h(\pi)|\geq 2\,comp(\pi)-1$$
 therefore the  total
 contribution of  terms   such that $J\not=\emptyset$ 
  can be bounded by $Cq^{n/2-1}$
 for some constant $C$ depending only on $n$ and on the sequence
 $(c_l)_{1\leq l\leq n}$. A similar argument using Lemma 4.3.2 of \cite{B}
  would apply to (2.1.3), we
 shall therefore restrict the sums in (2.1.2) and (2.1.3)
 to $i's$ and $j's$ in the range
 $1,\ldots, q$, and replace (2.1.4) by a sum
  $$
 \sum_{\pi}{1\over (q+1)^2}(q)_{comp(\pi)}\psi(\pi)\tag 2.1.6
 $$
 over partitions $\pi$ of $ \{1,2,\ldots,n_1\}\cup\{1',2',\ldots,n_2'\}$.
 Consider now the contribution to (2.1.6) of
 partitions $\pi$ such that some $u\in\{1,2,\ldots,n_1\}$ and some 
 $v'\in\{1',2',\ldots,n_2'\}$ are in the same component,
 then using Lemma 5.1.2 from \cite{B} and the estimate $(1)$ of 
 Definition 1.3 we can again
  bound this
 contribution by $Cq^{n/2-1}$.
 It remains to consider the contribution of partitions $\pi$ which split as the
 union of a partition $\pi_1$
  of $\{1,\ldots,n_1\}$ and a partition $\pi_2$ of $\{1',2',\ldots,
 n_2'\}$. For such a partition, any  associated permutation 
 $(i_1i_2)\ldots (i_{n_1}i_1)
 (j_1j_2)\ldots (j_{n_2}j_1)$ is the product of two
 permutations with disjoint supports 
 $(i_1i_2)\ldots (i_{n_1}i_1)$ and 
 $(j_1j_2)\ldots (j_{n_2}j_1)$, so that by the asymptotic
 factorization property, one has 
 $$\gather
 |\psi\bigl((i_1i_2)\ldots (i_{n_1}i_1)
 (j_1j_2)\ldots (j_{n_2}j_1)\bigr)-
 \psi\bigl((i_1i_2)\ldots(i_{n_1}i_1)\bigr)\psi\bigl(
 (j_1j_2)\ldots (j_{n_2}j_1)\bigr)|\\
 \leq \delta q^{-|h(\pi)|/2}
 \endgather
 $$
 Since by \cite{B}, Section 4.3,
  one has $n_i+|h(\pi_i)|\geq 2\,comp(\pi_i)-2$
  we can replace (2.1.6) by
 $$
 \sum_{\pi=\pi_1\cup\pi_2}{1\over (q+1)^2}
 (q)_{comp(\pi)}\psi(\pi_1)\psi(\pi_2)\tag 2.1.7
 $$
 making an error bounded by $C\delta q^{l/2}$.
 On the other hand the right hand side of (2.1.5) can be replaced by
 $$
 \sum_{(\pi_1,\pi_2)}{1\over (q+1)^2}(q)_{comp(\pi_1)}(q)_{comp(\pi_2)}
 \psi(\pi_1)\psi(\pi_2)\tag 2.1.8
 $$
The difference between (2.1.7) and (2.1.8) is 
$$
 \sum_{(\pi_1,\pi_2)}{1\over (q+1)^2}((q)_{comp(\pi)}-
 (q)_{comp(\pi_1)}(q)_{comp(\pi_2)})
 \psi(\pi_1)\psi(\pi_2)
 $$
  One has 
 $$comp(\pi_1\cup \pi_2)=comp(\pi_1)+comp(\pi_2)$$  thus 
 $$|(q)_{comp(\pi_1\cup\pi_2)}-
 (q)_{comp(\pi_1)}(q)_{comp(\pi_2)}|\leq Cq^{comp(\pi_1)+comp(\pi_2)-1}$$
 By Lemma
  4.3.2 and 4.3.3 of \cite{B}, one has 
  $$n_1+n_2+|h(\pi_1)|+|h(\pi_2)|\geq 2(comp(\pi_1)+comp(\pi_2))-4$$
  therefore 
  each term in this difference is bounded by
  $Cq^{l/2-1}$ and  one gets the
 stated result. \qed
 \enddemo
 Using the Lemma, we can now estimate the variance of $m_n(\lambda)$ under the
 measure $\Pi_{\psi}$, indeed it follows from Section 3.3 in \cite{B} that  one has
 $$
 \Pi_{\psi}((m_l(\lambda)-m_l(\psi))^2)=\tau_{\psi}(\langle
 \Gamma(\psi)^l\rangle^2)-\tau_{\psi}(\langle \Gamma(\psi)^l\rangle)^2$$
 therefore, using Lemma 2.1, for
 $\psi\in \Cal F_{c,\delta}^n$ one has
$$
 \Pi_{\psi}((m_l(\lambda)-m_l(\psi))^2) 
 \leq K(\delta+q^{-1})q^{l}\quad\text{by Lemma 2.1}$$
 Theorem 1 now follows by an application of Markov's inequality.\qed 
 \par
 Corollary 1 is a simple consequence of Theorem 1 and of the moment-cumulant
 formula. 
 \par
 Finally for the proof of Theorem 2, we can use the same arguments as
 in Section 4 of \cite{B}, and the hypotheses (1) and (2) of Definition 3.1.
 This shows that there exists a constant $K$ such that,
  for all $\psi\in\Cal F_{c,\delta}^n$, all $l\leq n$ and all cycle $\rho$
   on $l-1$
  elements,  one has
 $$|c_l(\psi)q^{1-l}
 -\psi(\rho)|\leq K(\delta+q^{-1}) q^{-l/2}\qquad\text{for all
 $l\leq n$}$$
 This, combined with Corollary 1 gives Theorem 2.\qed

 \head{3. Asymptotics of the tensor representation}\endhead
 \subhead{3.1 The case of the canonical trace}\endsubhead
 In this section we consider the action of $S_q$ on $(\Bbb
 C^N)^{\otimes q}$ by permutation of the factors in the tensor product. This
 representation plays a key role in the treatment by Schur and Weyl of the
 representation theory of both symmetric and general linear groups \cite{W}.
  We
 consider the decomposition of this representation into isotypic components
 $$(\Bbb
 C^N)^{\otimes q}=\bigoplus_{\lambda\in \Cal Y_q} E_{\lambda}$$
 then the
  relative dimensions $\dim E_{\lambda}\over N^q$ define a probability
 measure on $\Cal Y_q$. Let $\Pi_q$ be the image of this probability measure on 
 $\Cal Y_q\subset\Cal C\Cal Y$
 by the scaling map $\omega\mapsto q^{-1/2}\omega(q^{1/2}.)$. Thus $\Pi_q$ is
 a probability measure on the set $\Cal C\Cal Y$ of continuous Young
 diagrams.\par
 Let $\arcsin$ take values in $[-\pi/2,+\pi/2]$ and $\arccos$ in $[0,\pi]$,
 and define
 $$h(c,u)={2\over \pi}\biggl(u\arcsin\left[{u+c\over 2\sqrt{1+uc}}\right]
    +{1\over c}\arccos\left[{2+uc-c^2\over 2\sqrt{1+uc}}\right]
    +{1\over 2}\sqrt{4-(u-c)^2}\biggr)
    $$
    for $0<c<\infty$, and $u\in[c-2,c+2]$.
 Let us denote
 by $P_c$, the continuous diagram given by
 the formula
 $$\matrix P_0(u)=&\left\{\matrix
    {2\over \pi}(u\arcsin({u\over 2})+
    \sqrt{4-u^2})
    &\text{if $u\in [-2,+2]$}\\
    |u|&\text{if $u\notin[-2,+2]$}
 \endmatrix
\right.\hfill\,
 \\ \noalign{\vskip 5mm}
 \\
 P_c(u)=&
 \left\{
 \matrix
    h(c,u)
    &\text{if $u\in [c-2,c+2]$}\\
    |u|&\text{if $u\notin[c-2,c+2]$}
\endmatrix \right.\hfill\,
\\  \noalign{\vskip 2mm}
\\  &\text{for $0<c<1$}\hfill\,
\\  \noalign{\vskip 2mm}
\\
P_1(u)=&\left\{\matrix {u+1\over 2}+
{1\over \pi}\left((u-1)\arcsin({u-1\over
2})+\sqrt{4-(u-1)^2}\right)
&\text{if $u\in[-1,3]$}\\
|u| &\text{if $u\notin[-1,3]$}\endmatrix\right.
\\ 
 \noalign{\vskip 5mm}
\\
P_c(u)=&
 \left\{
 \matrix u+{2\over c}&\text{if $u\in[-{1\over c},c-2]$}\\
 h(c,u)&
    \text{if $u\in [c-2,c+2]$}\\
    |u|&\text{if $u\notin[-{1\over c},c+2]$}
 \endmatrix
\right.\hfill\,\\  \noalign{\vskip 2mm}
\\
&\text{for $c>1$}\hfill\,
\endmatrix\tag 3.1.1
$$
The diagrams are depicted below for various values of $c$.
$$\epsfxsize=12cm\epsfysize=8cm\epsfbox{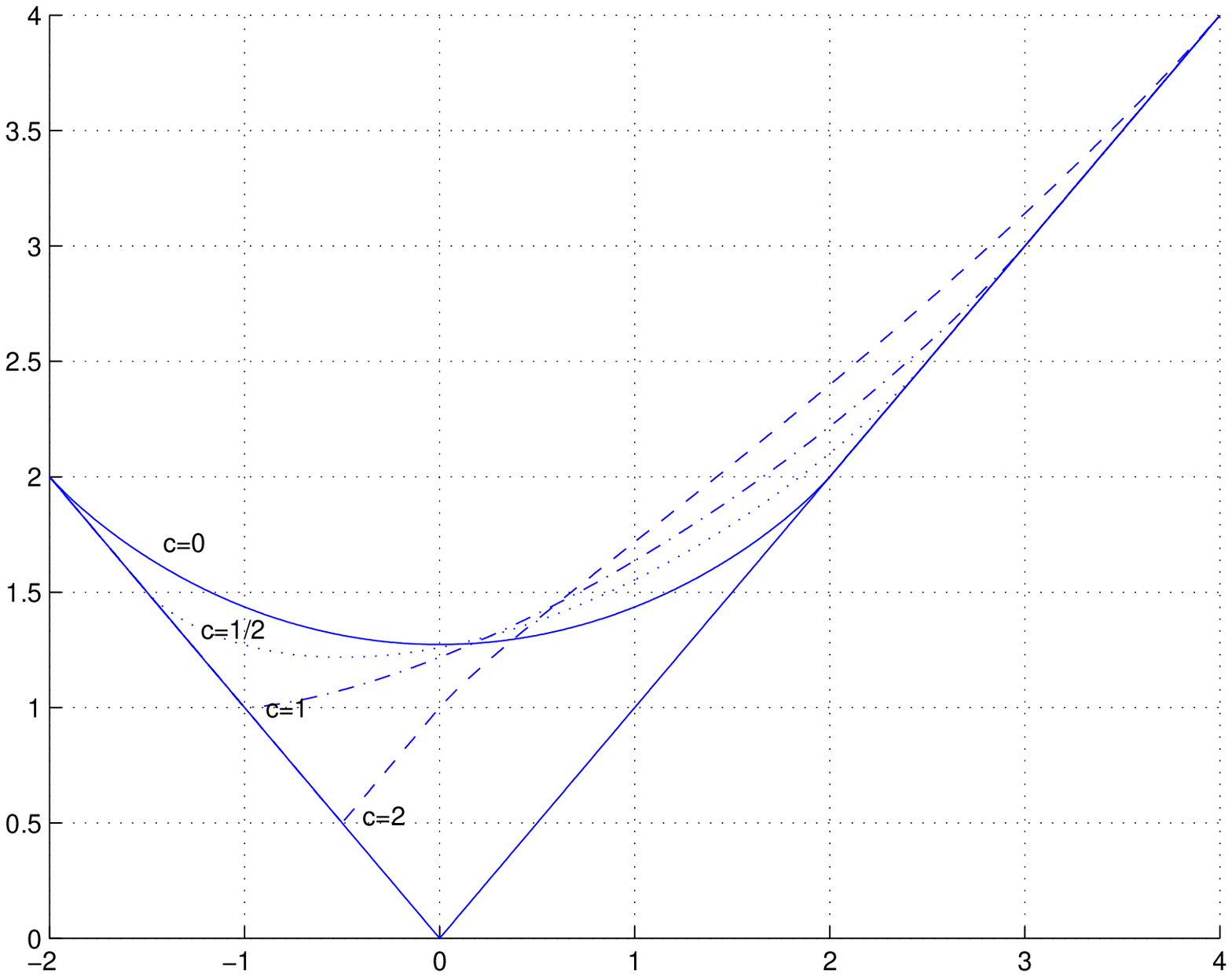}$$
The limit result for the $\Pi_q$ is the following. 

\proclaim{Theorem 3} 
As $q\to\infty$ and $\sqrt{q}/N\to c\in[0,+\infty[$, the measure $\Pi_q$
converges uniformly, in probability,
towards the Dirac measure  at the continuous diagram $P_c$.
\endproclaim
\demo{Proof}
The character of the tensor representation can be computed in the basis
$e_{i_1}\otimes \ldots \otimes e_{i_q}$, where $e_i$ is a basis of $\Bbb C^N$,
one gets, for $\rho\in S_q$,
$$tr(\rho)=N^{c(\rho)-q}=N^{-\vert\rho\vert}$$
This is a positive definite function on $S_q$. Taking $N=N(q)$ such that 
$\sqrt{q}/N\to c\in[0,+\infty[$, we see that the corresponding functions 
  belong to $\Cal F_{c,0}^n$
for all $n$, where the sequence $(c_l=C^l)_{l\geq
1}$,  with $C=\sup_q\{\sqrt{q}/N(q)\}$, 
therefore we can apply Theorem 1 and 2, and Corollary 1 
to see that under $\Pi_q$,  the sequence of  cumulants of a random diagram
converges in probability
 towards the sequence   
   $(0,1,c,c^2,\ldots, c^n,\ldots)$,
    and the sequence of moments converges towards
   the corresponding moments.  The correspondence between measures and diagrams
   given by (1.1) has the following continuity property (see \cite{K}): \par
   if 
  $\omega\in \Cal C\Cal Y$,
  then for all $\varepsilon>0$ there exists $n$ and $\delta$ such that, if
  $\omega'\in\Cal C\Cal Y$ and 
  $|m_l(\omega)-m_l(\omega')|<\delta$ for all $l\leq n$,  then
  $|\omega(u)-\omega'(u)\vert<\varepsilon$ for all $u\in\Bbb R$. 
  \par
  Therefore it is
  enough to check that the sequence $(0,1,c,c^2,\ldots, c^n,\ldots)$ is the
  sequence of cumulants of the diagram $P_c$. For $c=0$, $P_c$ is the diagram
  corresponding to the semi-circle distribution (see \cite{K}), while for 
    $c>0$
   the constant sequence 
  $(c^{-2},\ldots, c^{-2},\ldots)$ 
  is the sequence of cumulants of the so-called
  Pastur-Marcenko or free Poisson distribution (see \cite{PM}, or
  \cite{VDN}), with parameter $c^{-2}$. Taking the image of this measure by
  the mapping $x\mapsto cx-{1\over c}$  yields the measure with free cumulant
  sequence $(0,1,c,c^2,c^3\ldots)$. We shall now compute the
  corresponding diagram, according to \cite{K}. 
The $R$-transform 
with cumulant sequence $(0,1,c,c^2,c^3,\ldots)$ is 
$$R(z)=\sum_{n=1}^{\infty}z^nc^{n-1}={z\over 1-cz}$$
 thus 
  $$K(z)={1\over z}+R(z)=
  {1\over z}+{z\over 1-cz}$$ Inverting the function $K$ gives the Cauchy
 transform 
$$G(z)={z+c-\sqrt{(z-c)^2-4}\over 2(1+cz)}={2\over z+c+\sqrt{(z-c)^2-4}}$$
where the branch of the square root on $\Bbb C\setminus [0,+\infty[$
is chosen to have positive imaginary part. 
 We now compute the Rayleigh measure  $\tau$ 
according to the formula
$${\partial\over\partial z}\log G(z)=-\int_{\Bbb R} {\tau(du)\over z-u}$$
One has
$$
\aligned 
        -{\partial\over\partial z}\log G(z)
	&
	={1+{z-c\over \sqrt{(z-c)^2-4}}\over z+c+\sqrt{(z-c)^2-4}}
	\\
	&
	= {1\over 2\sqrt{(z-c)^2-4}}
	\left({2+cz-c^2+c\sqrt{(z-c)^2-4}\over 1+cz}\right)
\endaligned
$$
Let 
$$
k(c,u)={2+cu-c^2\over 2\pi(1+cu)\sqrt{4-(u-c)^2}}
\qquad\text{ 
for $-2+c<u<2+c$.}
$$
Using Stieltjes inversion formula we get
$$
\tau_c(du)=
     \left\{\matrix k(c,u)1_{[c-2,c+2]}(u)\,du
             &\text{ if $0\leq c<1$}\\
	     {1\over 2}\delta_{-1}(du)+
	     {du\over 2\pi\sqrt{(x+1)(3-x)}}1_{]-1,3[}(u)
	     &\text{ if $c=1$}\\
	     \delta_{-{1\over c}}(du)+k(c,u)1_{[c-2,c+2]}\,du&\text{ if $c>1$}
	     \endmatrix\right.\tag 3.1.2
	     $$
	     	  	     From the Rayleigh 
measure $\tau_c$ we can recover the diagram by
the formula 
$$P_c(u)=\int_{-\infty}^{+\infty}|u-x|\,\tau_c(dx)$$
A lenghty but straightforward computation using (3.1.2)  gives formula (3.1.1).
\qed
\enddemo
\subhead{3.2 Generalization to other tensor states}\endsubhead
 Let $H=H_0\oplus H_1$, be a  $\Bbb Z/2$-graded Hilbert space 
with grading operator $\Delta$.
Endow $\Cal B(H^{\otimes q})=\Cal B(H)^{\otimes q}$ with a tensor state
$\tau^{\otimes q}$ where $\tau$ is a state on $\Cal B(H)$, of the form
$\tau(X)=Tr(TX)$, with $T$ a positive  trace class operator, with trace 1,
commuting with the grading. 
Denote
by $t_{1,0}\geq t_{2,0}\geq \ldots\geq t_{k,0}\geq \ldots\geq 0$ and 
$t_{1,1}\geq t_{2,1}\geq \ldots\geq t_{k,1}\geq \ldots\geq 0$  
the eigenvalues of $T$ on $H_0$ and $H_1$ respectively, and
let $p_n=Tr(\Delta(\Delta T)^n)=\sum_{j=1}^{\infty}t_{j,0}^n-(-t_{j,1})^n$ for
$n\geq 1$.

In this last section we shall consider the representation of $S_q$ on
$H^{\otimes q}$, given by 
$$(i\,i+1)(v_1\otimes \ldots v_i\otimes v_{i+1} \ldots\otimes v_q)
=(-1)^{\varepsilon(v_i)\varepsilon
(v_{i+1})}v_1\otimes\ldots v_{i+1}\otimes v_{i} \ldots\otimes v_q$$
where $v_k$ are graded vectors, with graduation $\varepsilon(v_k)\in\{0,1\}$. 
The positive
definite function on $S_q$ determined by this representation and the
 state $\tau^{\otimes q}$ is given
by
$$\tau^{\otimes q}(\rho)=\prod_{c|\rho}p_{\vert c\vert+1}\tag 3.2.1$$
where the product is over the non trivial cycles $c$ of $\rho$. Recall that for
a cycle $c$ one has $|c|+1=s(c)$.
This is of
course a direct generalization of Section 3.1 where we had $t_{1,0}
=t_{2,0}=\ldots
=t_{N,0}={1\over N}$, the other eigenvalues being 0. We shall now investigate the
asymptotic decomposition of the associated probability measure on Young
diagrams, under the hypothesis that  $\tau:=\tau_q$ depends on $q$ as
$q\to\infty$, and that the moments $Tr(T_q^n)$ satisfy
$$q^{(n-1)/2}Tr(T_q^n) \to_{q\to\infty} w_{n+1}\quad n\geq 1$$
for some sequence of real numbers $w_n$, with $|w_n|\leq C^{n}$ for some
constant $C>0$. This last condition insures that the positive definite functions
(3.2.1) belong to the sets $\Cal F^n_{c,0}$ as in Section 3.1.
  Using
the same arguments as in the preceding section, we see that the measure
on rescaled diagrams converges towards the Dirac mass at the diagram with
cumulants
$(w_n)_{n\geq 1}$, with $w_1=0,\,w_2=1$.
 These cumulants have the form $w_n=\int_{\Bbb R} |t|t^{n-2}\mu(dt);\ n\geq 2$
 for some positive measure $\mu$ on $[-C,C]$,
satisfying $\int_{\Bbb R} |t|\mu(dt)=1$. The $R$-transform of the measure
with cumulants $w_n$ is thus 
$$R(z)=\sum_{n=1}^{\infty} w_nz^{n-1}=\int_{\Bbb R}{z\over 1-zt}|t|\mu(dt)$$
Comparing with the free L\'evy-Khintchine formula of
\cite{BV}, we see that the  measure corresponding to the limit diagram
is freely
indefinitely divisible, with corresponding L\'evy measure $\nu(dt)={1\over |t|}
\mu(dt)$
\proclaim{Remark} One can give an explicit expression for the weight of a given
diagram $\lambda\in\Cal Y_q$, under the measure associated with the function
(3.2.1), in terms of
the Schur functions (see also \cite {M} formula (7.8)). It would be interesting
to rederive the results of Sections 3.1 and 3.2 directly from this explicit
expression.\endproclaim

\refstyle{A}


\widestnumber \key{VDN}
\Refs\nofrills{References}

\ref\key B\by P. Biane\paper Representations of symmetric groups and free
probability\jour Adv. Math.\yr 1998\pages 126--181
\vol 138\endref
\ref \key BV
\by  H. Bercovici and D. Voiculescu
\paper  Free convolution of
measures with unbounded support
\jour Indiana University Mathematics
Journal
\vol  42 
\yr 1993 
\pages  733-773\endref
\ref \key K\by S. V. Kerov\paper Interlacing measures
\yr 1998\inbook  Kirillov's seminar on representation theory
\bookinfo Amer. Math. Soc. Transl. Ser. 2
Appl.\pages 35--83\vol 181\publ Amer. Math. Soc.\publaddr
Providence, RI\endref
\ref \key LS\by B. F. Logan and L. A. Shepp
\paper A variational problem for random Young tableaux
\jour Adv.  Math.\vol 26\yr 1977\pages 206--222\endref
\ref\key Mac\by I. G. Macdonald\book Symmetric functions and Hall
polynomials, $2^{nd}$ Ed.\publ Oxford Univ. Press\publaddr Oxford\yr 1995 \endref
\ref\key MP\by V.A.
Mar\v cenko and  L. A. Pastur\paper 
Distribution of eigenvalues in certain sets of random matrices \jour
Mat. Sb. (N.S.) \vol 72 (114)\yr 1967\pages 507--536\endref 
\ref\key S\by R. Speicher\book Combinatorial theory of the free product with
amalgamation and operator-valued free probability theory\bookinfo Mem.
A.M.S.\vol 627
\publ Amer. Math. Soc.\publaddr
Providence, RI\yr 1998\endref
\ref \key VK\by A. M. Vershik and S. V. Kerov\paper Asymptotics of the
Plancherel measure of the symmetric group \jour Soviet Math. Dokl.  \vol 18 
\yr 1977\pages 527--531\endref
\ref\key VDN\by  D. V. Voiculescu,  K. Dykema
and  A. Nica\book Free random variables
\bookinfo  CRM Monograph Series No. 1\publ Amer. Math. Soc.\publaddr
Providence, RI\yr 1992\endref
\ref\key Wey \by H. Weyl\book The Classical Groups. Their invariants and
Representations\publ Princeton Univ. Press\yr 1939
\publaddr Princeton N.J.\endref
\endRefs
\enddocument